\documentclass[12pt,a4paper,twoside]{article}

\usepackage[centertags]{amsmath}
\usepackage[T1]{fontenc}
\usepackage{amsfonts}
\usepackage{amssymb}
\usepackage{ntheorem}
\usepackage{ulem}
\usepackage{epsfig}
\usepackage{subfigure}
\usepackage[all]{xy}
\usepackage{changebar}
\usepackage{color}
\usepackage{mathrsfs}
\usepackage[boldsans]{ccfonts}
\usepackage{array}
\usepackage{tabulary}
\usepackage{supertabular}
\usepackage{calc}
\usepackage{eepic}
\usepackage{titlesec}
\usepackage{graphicx}
\usepackage{enumitem}

\theoremstyle{break} \newtheorem{theorem}{Theorem}[section]
\theoremstyle{break} 
\theoremstyle{break} 
\theoremstyle{nonumberbreak}  
\theoremstyle{break} \newtheorem{lemma}[theorem]{Lemma}
\theoremstyle{break} 
\theoremstyle{break} 
\theoremstyle{break} 
{\theorembodyfont{\rmfamily}
\newtheorem{remark}[theorem]{Remark}
{\theorembodyfont{\rmfamily}}
\theoremstyle{break} 
\theoremstyle{break} 
\theoremstyle{break} 
\theoremstyle{break} 
\numberwithin{equation}{section}

\newcommand{\hide}[1]{}

\newcommand{\R}{{\mathbb{R}}}
\newcommand{\D}{{\mathbb{D}}}
\newcommand{\C}{{\mathbb{C}}}

\def\Re{\mathop{{\rm Re}}}

\def\Aut{\mathop{{\rm Aut}}}

\begin{document}

\begin{center}
{\Large \bf Composition and decomposition of indestructible\\[2mm]  Blaschke products}
\footnotetext{2000 Mathematics Subject
 Classification: 30H05; 30J05, 30J10.\\ This research has been
supported by the Deutsche Forschungsgemeinschaft (Grants: Ro 3462/3--1 and Ro
3462/3--2) \hfill{Date: \today}}

\end{center}
\renewcommand{\thefootnote}{\arabic{footnote}}
\setcounter{footnote}{0}

\smallskip

\begin{center}
{\large Daniela Kraus and Oliver Roth}\\[3mm]
\end{center}

\smallskip

\begin{center}
\begin{minipage}{13.2cm}
{\bf Abstract.}{ \small
We prove that the composition of two indestructible Blaschke products is again
an indestructible Blaschke product. We also show that if an indestructible
Blaschke product is the composition of two bounded analytic functions, then
both functions are indestructible Blaschke products.}
 \end{minipage}
\end{center}

\medskip

\section{Introduction and results}

Let $H^{\infty}$ denote the Banach space of all functions analytic and bounded
in the 
unit disk  $\D:=\{z \in \C \, : \, |z|<1\}$ equipped with the supremum norm
$||\cdot||_{\infty}$. The set $\mathcal{B}:=\{f \in H^{\infty} \text{ nonconstant}\, : \,
||f||_{\infty}  \le 1\}$ %(modulo the functions in $H^{\infty}$ of constant
%absolute value $1$) 
of nonconstant functions in the unit ball of $H^{\infty}$
is clearly closed under composition, that is,
$$ f \circ g \in \mathcal{B} \qquad \text{ for all } \quad f,g \in \mathcal{B} \,. $$
There are three well--known subsemigroups inside the composition semigroup $\mathcal{B}$:
\begin{itemize}
\item[(a)] the set
$$\Aut(\D)=\left\{\eta \frac{z-a}{1-\overline{a} z} \, :\,  |a|<1 \, ,
|\eta|=1 \right\}$$
of all conformal automorphisms of $\D$. Note that $\Aut(\D)$ is actually a group w.r.t.~composition.
\item[(b)] the set of all \textit{finite} Blaschke products
$$   \eta \prod \limits_{j=1}^N \frac{z-a_j}{1-\overline{a}_j z}\, ,
  \qquad  a_1,\ldots, a_N \in \D, \, |\eta|=1 \, ;$$
\end{itemize}
and
\begin{itemize}
\item[(c)] the set of all inner functions, i.e., those
  functions $F \in H^{\infty}$ for which the radial limit function
$$ F^*(\zeta):=\lim \limits_{r \to 1-} F(r \zeta)\, , $$
%which exists for almost every (a.e.) $\zeta \in \partial \D$ by Fatou's
%theorem, 
satisfies $|F^*(\zeta)|=1$ for a.e.~$\zeta \in \partial \D$.
We refer to \cite{Ryf1966} for a proof of the fact that the set of inner functions 
is a semigroup and to \cite{GLMR1994} for more about the structure of this semigroup.
\end{itemize}

On the other hand, the set of all Blaschke products\footnote{We use the
  convention $\frac{a}{|a|}=1$ if $a=0$.}
$$ B(z)=\eta  \prod \limits_{j=1}^{\infty} \frac{-\overline{a_j}}{|a_j|}
\frac{z-a_j}{1-\overline{a}_j z} \, , \quad |\eta|=1\, , \quad (a_j)
\subseteq \D \text{ s.t. } \sum \limits_{j=1}^{\infty} \left(1-|a_j|
\right)<\infty \ ,$$
is \textit{not} closed under composition. In fact, a celebrated result of
Frostman \cite[Theorem 6.4]{Gar2007} says that every inner function can be
written as the
composition $T \circ B$ of a unit disk automorphism $T$ and a Blaschke
product $B$. In view of Frostman's result it is natural to 
consider those Blaschke products $B$ such that $T \circ B$ is a Blaschke product
for \textit{every} $T \in \Aut(\D)$. Such Blaschke products are called 
\textit{indestructible}, see \cite{McL1972,Ros2008}.
Clearly, every finite Blaschke product is indestructible, but there are also
numerous examples of infinite indestructible Blaschke products (see
\cite[Chapter 5]{Ros2008} and Remark \ref{rem:maximal} below).
Indestructible Blaschke products do have a number of intriguing properties,
some of which are described e.g.~in \cite{McL1972, Ros2008}.

\medskip

Our first result shows that the set
of all indestructible Blaschke products forms a composition semigroup.

\begin{theorem}[Composition of indestructible Blaschke products] \label{thm:1}
The composition of two indestructible Blaschke products is an indestructible
Blaschke product.
\end{theorem}

Knowing that $B \circ C$ is an indestructible Blaschke product if both
$B$ and $C$ are indestructible, the second result of this paper deals
with the ``inverse'' problem.

\begin{theorem}[Decomposition of indestructible Blaschke products] \label{thm:2}
Let $B,C \in H^{\infty}$ be nonconstant functions such that $B \circ C$ is an indestructible Blaschke
product. Then $B$ and $C$ are both indestructible Blaschke products.
\end{theorem}

The proofs of Theorem \ref{thm:1} and Theorem \ref{thm:2} are given in the
next two sections.

\begin{remark}[Prime indestructible Blaschke products]
In view of Theorem \ref{thm:1} and Theorem \ref{thm:2}, it is natural to
ask for \textit{prime} indestructible Blaschke products, i.e.~indestructible
Blaschke products $A$ such that if $A=B \circ C$ for $B,C \in H^{\infty}$,
then $B$ or $C$ is a unit disk automorphism. It would be interesting to find
examples of infinite prime indestructible Blaschke products. We note that
the concept of primality in the context of $H^{\infty}$--functions 
has been considered for a long time
for finite Blaschke products (see \cite{NgTsang2012,NgWang2011} for recent developments).
For inner functions the concept is due to Stephenson \cite{Ste1982} (see also
\cite{Ste1977,Ste1978}) and has been further explored 
e.g.~by Gorkin, Laroco, Mortini and Rupp \cite{GLMR1994}.
\end{remark}

\begin{remark}[Maximal Blaschke products] \label{rem:maximal}
A result intimately connected to Theorem \ref{thm:1} has recently
been proved  in \cite{KR2012}, see also \cite{KR2012a}. There,  so--called maximal
Blaschke products have been studied. Maximal Blaschke products are
characterized by an extremal property and constitute an appropriate generalization of the class of finite
Blaschke products. They are defined as follows. Let $F \in
H^{\infty}$ be a  
nonconstant function, let $\mathcal{C}:=(z_1,z_2,\ldots)$ denote the critical points of
$F$ counting multiplicity, and let $N$ denote the number of times that $0$
appears in $\mathcal{C}$. Consider the extremal problem
$$ \max \left\{\Re f^{(N+1)}(0) \, : \, f \in H^{\infty}, \, f'(z)=0 \text{
    for any } z \in \mathcal{C}\right\} \, .$$
It is not difficult to see that this extremal problem
has a unique solution $F_{\mathcal{C}} \in H^{\infty}$ and that $F'_{\mathcal{C}}$
vanishes precisely on the sequence $\mathcal{C}$.   
It turns out, see \cite[Theorem 1.1]{KR2012}, that the extremal function
 $F_\mathcal{C}$ is an indestructible Blaschke product. This generalizes Nehari's
well--known extension \cite{Neh1946} of the Ahlfors--Schwarz lemma from 1947.
Nehari's result covers \textit{finite} sequences $\mathcal{C}$, in which cases the extremal
functions are precisely the \textit{finite} Blaschke products.  Now, 
every Blaschke product of the form $T \circ F_{\mathcal{C}}$ with $T \in
\Aut(\D)$ and $F_{\mathcal{C}}$ the extremal function for the critical set
$\mathcal{C}$ of some nonconstant $H^{\infty}$--function  is called a
\textit{maximal Blaschke product}. In particular, every finite Blaschke
product is maximal  and every maximal Blaschke product is  
indestructible, so maximal Blaschke products provide a
large collection of examples of indestructible Blaschke products.
Maximal Blaschke products do have a number
of striking properties reminiscent of finite Blaschke products and
Bergman space inner functions and they are intimately connected to hyperbolic
geometry, see  \cite{Kra2013,KR2012a,KR2012}. 
From the point of view of the present paper, the perhaps most interesting property is 
the fact that  the set of maximal Blaschke products is closed under
composition, see \cite[Theorem 1.7]{KR2012} and \cite[Theorem 7]{KR2012a}.
\end{remark}

\begin{remark}
Summarizing the results and remarks above,
 we have the following ``zoo'' of subsemigroups of 
$\mathcal{B} \subset H^{\infty}$ involving Blaschke products (BPs):

\medskip
\fbox{\parbox[h][\height -11mm +\baselineskip][c]{15.3cm}{
$$ \Aut(\D) \subset \{\text{finite BPs}\} \subset 
\{\text{max.~BPs}\} \subset \{\text{indestr.~BPs}\} \subset
\{\textit{inner functions\}} \subset \mathcal{B} \, .
$$}}

\bigskip

All but possibly one of these inclusions are strict: We do not know
whether there exists an indestructible Blaschke product which is not a
maximal Blaschke product.
\end{remark}

Pioneering work on indestructible Blaschke products can be found
in the papers by Heins \cite{Hei1953,Hei1955}, McLaughlin \cite{McL1972} and
Morse \cite{Mor1980}.
We also refer to the excellent survey on indestructible Blaschke products by Ross
\cite{Ros2008} and to the Fields Institute Proceedings \cite{FM2012} for a
collection
of surveys and research articles on
Blaschke products and inner functions in general.

\section{Proof of Theorem \ref{thm:1}} 
The proof of Theorem \ref{thm:1} relies on the following characterization of
indestructible Blaschke products due to McLaughlin \cite[Theorem 1]{McL1972}; see also
Ross \cite[Theorem 4.1]{Ros2008}. We use the following notation. If $F \in
H^{\infty}$ is not constant, then for fixed $a \in \D$ let $\xi_1(F;a),\xi_2(F;a),\ldots \in \D$ denote
the solutions to $F(z)=a$ counting multiplicities.
We denote by  $z_1(F),z_2(F),\ldots \in
\D\backslash \{0\}$ 
the \textit{non--zero} solutions to $F(z)=F(0)$ again counting multiplicities.

\begin{lemma} \label{lem:ml}
Let $F \in H^{\infty}$  be an inner function such that 
$$ F(z)=F(0)+b_n z^n+b_{n+1} z^{n+1}+\ldots \qquad (b_n\not=0) \, . $$ 
 Then  $F$ is an indestructible Blaschke product if and only if the following
two conditions hold:
\begin{eqnarray} \label{eq:m1}
 \left|\frac{F(0)-a}{1-\overline{a} \, F(0)}\right| &=& \prod \limits_{j=1}^{\infty} |\xi_j(F;a)| \quad \text{ for any }
a \in \D \backslash \{ F(0)\} \\ \label{eq:m2}
\frac{|b_n|}{1-|F(0)|^2} &=& \prod \limits_{j=1}^{\infty} |z_j(F)| \, . 
\end{eqnarray}
\end{lemma}

\begin{remark}
 The statement of Lemma
  \ref{lem:ml} in \cite{McL1972} and \cite{Ros2008} is slightly weaker, since 
both sources assume that $F$ is a Blaschke product. In fact, 
the proofs in \cite{McL1972,Ros2008} reveal that it is sufficient to assume that
$F$ is an inner function. This simple observation will be crucial in the proof of
Theorem \ref{thm:1} given below.
\end{remark}

{\it Proof of Theorem \ref{thm:1}.}
Let $B$ and $C$ be two indestructible Blaschke products. We can assume
 that both $B$ and $C$ are not constant.
%with
%$$ B(z)=B(0)+b_N z^N+b_{N+1} z^{N+1}+\ldots \, , \qquad C(z)=C(0)+c_M
%z^{M}+c_{M+1} z^{M+1}+ \ldots \, .$$
Since $B$ and $C$ are inner functions, we see that $A:=B \circ C$ is an inner function.
Fix $a \in \D$.
Since $B$ is indestructible, the function $$w \mapsto \frac{B(w)-a}{1-\overline{a} B(w)}$$ is a
Blaschke product. Recalling the convention $\frac{\xi}{|\xi|}=1$ if $\xi=0$,
this implies
\begin{equation} \label{eq:0}
 \frac{B(w)-a}{1-\overline{a} B(w)} =\eta_a
 \prod \limits_{j=1}^{\infty} \frac{-\overline{\xi_j(B;a)}}{|\xi_j(B;a)|}
\, \frac{w-\xi_j(B;a)}{1-\overline{\xi_j(B;a)} w} \, , \qquad w \in \D \, , 
\end{equation}
for some $\eta_a \in \partial \D$.
We now  distinguish the two cases 
\begin{center}
\begin{tabular}{rl}
Case I: &  $a\not=A(0)$ \\[1mm]
Case II: & $a=A(0)$.
\end{tabular}
\end{center}
\medskip

We start with Case I, so  $a \in \D \backslash \{A(0)\}$.
For $w=C(0)$ equation (\ref{eq:0}) gives
\begin{equation} \label{eq:1}
 \frac{A(0)-a}{1-\overline{a} A(0)} =
\displaystyle \eta_a \prod \limits_{j=1}^{\infty} \frac{-\overline{\xi_j(B;a)}}{|\xi_j(B;a)|}
\, \frac{C(0)-\xi_j(B;a)}{1-\overline{\xi_j(B;a)} C(0)} 
\end{equation}
As $A(0)\not=a$, this implies in particular,
$$ C(0)\not\in \{\xi_1(B;a),\xi_2(B;a),\ldots\}  \, .
$$
Since  $C$ is indestructible, Lemma \ref{lem:ml} therefore shows that
$$  \prod \limits_{k=1}^{\infty} |\xi_k(C;\xi_j(B;a))| =
\displaystyle \left| \frac{C(0)-\xi_j(B;a)}{1-\overline{\xi_j(B;a)} C(0)}
\right| \, ,\qquad j=1,2,\ldots \, , $$
so (\ref{eq:1}) takes the form
\begin{equation} \label{eq:2} \left| \frac{A(0)-a}{1-\overline{a} A(0)} \right|= \prod
\limits_{j,k=1}^{\infty} |\xi_k(C;\xi_j(B;a))|\, . 
\end{equation}
We now consider the equation $A(z)=a$. Let $z \in \D$ such that $A(z)=a$. This
is equivalent to $B(C(z))=a$, which is the same as 
$$ C(z) \in
 \{\xi_1(B;a),\xi_2(B;a),\ldots\} \, .
$$
Hence
$$ \{\xi_l(A;a) \, : \,
l=1,2\ldots\}= \{\xi_k(C;\xi_j(B;a)) \, : \,k,j =1,2,\ldots\} \, ,$$
so 
$$
 \prod \limits_{l=1}^{\infty} |\xi_l(A;a)|= \prod \limits_{j,k=1}^{\infty}
|\xi_k(C;\xi_j(B;a))| \, .
$$
Inserting this expression into (\ref{eq:2}) leads to 
$$ \left|
  \frac{A(0)-a}{1-\overline{a} A(0)} \right|=
\prod \limits_{l=1}^{\infty} |\xi_l(A;a)|\, .$$
Hence we have verified that condition (\ref{eq:m1}) holds for $F=A$.

\medskip

Let us turn to Case II, so assume $a=A(0)$. Now, we need to distinguish two subcases

\begin{center}\begin{tabular}{rl}
Case IIa: &  $C(0)\not=0$ \\[1mm]
Case IIb: & $C(0)=0$.
\end{tabular}\end{center}

\medskip

Let first $C(0)\not=0$ and suppose 
\begin{equation*}
\begin{array}{rcll}
C(z)&= & C(0)+c_N z^N+c_{N+1} z^{N+1}+\ldots\,,  \qquad & c_N\not=0\,,\\[2mm]
B(z) &=& B(0)+b_M z^M+b_{m+1} z^{M+1}+\ldots\,, \qquad  & b_M\not=0 \, .
\end{array}
\end{equation*}
Then $A(z)=B(C(z))$ has an expansion of the form
$$A(z)=A(0)+a_N z^N+a_{N+1} z^{N+1}+\ldots\,  $$
about $z=0$ with $$a_N=Nc_N b_M C(0)^{N-1}\not=0\, .$$
Since $a=A(0)=B(C(0))$, we have  
$$C(0) \in \{\xi_j(B;a)\, : \, j=1,2\ldots\}\, , $$ so
we may assume that $C(0)=\xi_1(B;a)$.
Using (\ref{eq:0}) for $w=A(z)$ with  $z\in \D \backslash\{0\}$
and dividing by $z^N$, we get
\begin{equation*} \label{eq:3}
\frac{1}{z^N} \frac{A(z)-a}{1-\overline{a} A(z)}=\eta_a \prod \limits_{j=2}^{\infty} \frac{-\overline{\xi_j(B;a)}}{|\xi_j(B;a)|}
\, \frac{C(z)-\xi_j(B;a)}{1-\overline{\xi_j(B;a)} C(z)} 
\cdot \frac{-\overline{C(0)}}{|C(0)|}
\, \frac{C(z)-C(0)}{1-\overline{C(0)} C(z)} 
\frac{1}{z^N}\, .
\end{equation*}
Letting $z \to 0$ on both sides, we obtain
\begin{equation} \label{eq:4}
 \frac{a_N}{1-|a|^2} =\eta_a \prod \limits_{j=2}^{\infty} \frac{-\overline{\xi_j(B;a)}}{|\xi_j(B;a)|}
\, \frac{C(0)-\xi_j(B;a)}{1-\overline{\xi_j(B;a)} C(0)} \cdot
\frac{c_N}{1-|C(0)|^2} \, . 
\end{equation}
In particular, since $a_N\not=0$ and $c_N \not=0$, we have
\begin{equation} \label{eq:5}
C(0)\not\in \{\xi_j(B;a) \, : \, j=2,3,\ldots\} \, .
\end{equation}
We now use the assumption that $C$ is indestructible. Lemma \ref{lem:ml} shows that
$$ \left| \frac{C(0)-\xi_j(B;a)}{1-\overline{\xi_j(B;a)} C(0)} \right|=\prod
\limits_{k=1}^{\infty} |\xi_k(C;\xi_j(B;a))| \, , \qquad j=2,3, \ldots \, , $$
and 
$$ \frac{|c_N|}{1-|C(0)|^2}=\prod \limits_{j=1}^{\infty} |z_j(C)| \, .$$
Inserting these last two expressions into (\ref{eq:4}), we arrive at
\begin{equation} \label{eq:6}
 \frac{|a_N|}{1-|a|^2}= 
\prod \limits_{j=2}^{\infty} \prod
\limits_{k=1}^{\infty} |\xi_k(C;\xi_j(B;a))| \cdot \prod \limits_{j=1}^{\infty} |z_j(C)|
\, .
\end{equation}

Now, let us find the \textit{non--zero} solutions to $A(z)=a$, i.e., the points
$z_j(A)$, $j=1,2\ldots$.
If $z \in \D\backslash \{0\}$ with $A(z)=a$, then $B(C(z))=a$. If $C(z)=C(0)$, then $z=z_j(C)$
for some $j=1,2\ldots$ and any $z_j(C)$ is a solution to $A(z)=a$.
If $C(z)\not=C(0)$, then (\ref{eq:5}) implies 
$z=\xi_k(C;\xi_j(B;a))$ for some $k=1,2,\ldots$ and
some $j=2,3,\ldots$. Conversely, each such $\xi_k(C;\xi_j(B;a))$ is a solution to
$A(z)=a$. Hence, we have shown that
\begin{eqnarray*} \{z_j(A) \, : \, j=1,2,\ldots\} &=& \\
&& \hspace*{-3cm} = \{z_j(C) \, : \, j=1,2,\ldots\} \cup  \{ 
\xi_k(C;\xi_j(B;a)) \, : \, k=1,2, \ldots, j=2,3, \ldots\}
\, . \end{eqnarray*}
This enables us to rewrite (\ref{eq:6}) as 
$$ \frac{|a_N|}{1-|a|^2}= \prod \limits_{j=1}^{\infty} |z_j(A)| \, .$$
Hence condition (\ref{eq:m2}) holds for $F=A$ in the case $C(0)\not=0$.

\bigskip

In a final step, we now proceed to establish condition (\ref{eq:m2})  for $F=A$  in
the remaining case $C(0)=0$. Let 
\begin{equation*}
\begin{array}{rcll}
B(z)&=&B(0)+b_N z^N+b_{N+1}z^{N+1}+\ldots \, , \hspace*{1cm} & b_N\not=0, \, N \ge
1\, , \\[2mm]
C(z)&=&c_M z^M+c_{M+1} z^{M+1}+\ldots\, , & c_M\not=0, \, M \ge 1\,.
\end{array}
\end{equation*}
Then $A(0)=B(C(0))=B(0)=a$ and
$$ A(z)=a+a_{NM}  z^{NM}+\ldots \, , \qquad a_{NM}=b_N c_M\not=0 \, . $$
Since $B$ is indestructible and $B(z)-a$ has a zero of order $N$ at $z=0$, we
get
$$ \frac{B(w)-a}{1-\overline{a} B(w)}=\eta_a w^N \prod \limits_{j=1}^{\infty} 
\frac{-\overline{z_j(B)}}{|z_j(B)|} \frac{w-z_j(B)}{1-\overline{z_j(B)} w} \,
, \qquad w \in \D \, .$$
For $w=C(z)$ and $z\in \D \backslash \{0\}$, this leads to
$$ \frac{1}{z^{MN}} \frac{A(z)-a}{1-\overline{a} A(z)} =\eta_a
\left(\frac{C(z)}{z^M} \right)^N \prod \limits_{j=1}^{\infty} 
\frac{-\overline{z_j(B)}}{|z_j(B)|} \frac{C(z)-z_j(B)}{1-\overline{z_j(B)} C(z)} \,
. $$
Letting $z \to 0$, we deduce
\begin{equation} \label{eq:7} 
 \frac{|a_{NM}|}{1-|a|^2}= |c_M|^N \prod \limits_{j=1}^{\infty} |z_j(B)| \, .
\end{equation}
Since $C$ is indestructible, Lemma \ref{lem:ml} implies
$$ |c_M|=\frac{|c_M|}{1-|C(0)|^2}=\prod \limits_{j=1}^{\infty} |z_j(C)| \, .$$
as well as 
$$ |z_j(B)|= \left| \frac{C(0)-z_j(B)}{1-\overline{z_j(B)} C(0)}
\right|=\prod \limits_{k=1}^{\infty} |\xi_k(C;z_j(B))| \, , \qquad j=1,2
\ldots \, .$$
Inserting the last two expressions into (\ref{eq:7}), we get
\begin{equation} \label{eq:8}
 \frac{|a_{NM}|}{1-|a|^2}= \prod \limits_{j=1}^{\infty} |z_j(C)|^N \cdot 
 \prod \limits_{j=1}^{\infty} \prod \limits_{k=1}^{\infty} |\xi_k(C;z_j(B))| \, .
\end{equation}
Consider the equation $A(z)=a$ and its \textit{non--zero} solutions.
Let $z\in \D \backslash \{0\}$ with $B(C(z))=A(z)=a$. If $C(z)=0$, then $z=z_j(C)$ is a zero of
$A(z)-a$ of order $N$. If $C(z)\not=0$, then $z=\xi_k(C;z_j(B))$ for
$j,k=1,2\ldots$. 
Hence we can write (\ref{eq:8}) as
$$ \frac{|a_{NM}|}{1-|a|^2}= \prod \limits_{j=1}^{\infty} |z_j(A)| \, .$$
This proves (\ref{eq:m2})  for $F=A$ also in
the case $C(0)=0$.

\medskip

In summary, we have shown that conditions (\ref{eq:m1}) and (\ref{eq:m2}) are
satisfied for $F=A$. Lemma \ref{lem:ml} therefore guarantees that $A$ is an
indestructible Blaschke product. The proof of Theorem \ref{thm:1} is complete.

\section{Proof of Theorem \ref{thm:2}}

We need the following well--known characterization of Blaschke products, see
\cite[Theorem 2.4]{Gar2007}.

\begin{lemma} \label{lem:2}
Let $f \in H^{\infty}$, $||f||_{\infty} \le 1$. Then the following are
equivalent.
\begin{itemize}
\item[(a)] $f$ is a Blaschke product.
\item[(b)] $\lim \limits_{r \to 1} \int \limits_{0}^{2 \pi} \log |f(r e^{i
    t})| \, dt=0$.
\item[(c)] The least harmonic majorant of $\log |f|$ is $0$.
\end{itemize}
\end{lemma}

{\it Proof of Theorem \ref{thm:2}.} 
We first prove that $C$ is an indestructible Blaschke product.
Let $T$ be a unit disk automorphism.  We need to  show that $\tilde{C}:=T \circ C$ is
a Blaschke product. Since $B$ is not constant, we can 
choose another unit disk automorphism $S$ such that
$\tilde{B}:=S \circ B \circ T^{-1}$ maps $0$ to $0$. Hence the Schwarz lemma
implies
\begin{equation} \label{eq:s}
|\tilde{B}(z)| \le |z| \, , \qquad z \in \D \, .
\end{equation}

Since $A:=B \circ C$ is indestructible, the function  $\tilde{A}:=S \circ A$ is a Blaschke
product and  $\tilde{A}=\tilde{B} \circ \tilde{C}$.
Using Lemma \ref{lem:2}, (a) $\Rightarrow$ (b), for $f=\tilde{A}$, we obtain
\begin{eqnarray*}
 0 = \lim \limits_{r \to 1} \int \limits_{0}^{2\pi} \log |\tilde{A}(r e^{it})|=\lim \limits_{r \to 1} \int \limits_{0}^{2\pi} \log |\tilde{B}(\tilde{C}(r e^{it}))|
 \, dt 
 \overset{(\ref{eq:s})}{\le} \lim \limits_{r \to 1} \int \limits_{0}^{2\pi} \log |\tilde{C}(r e^{it})|
 \, dt \le 0 \, .
\end{eqnarray*}
Hence
$$ \lim \limits_{r \to 1} \int \limits_{0}^{2\pi} \log |\tilde{C}(r e^{it})|
 \, dt =0 \, ,$$
so by Lemma \ref{lem:2}, (b) $\Rightarrow$ (a), we conclude that $\tilde{C}$ is a Blaschke product.

\medskip

Now, let us show that $B$ is an indestructible Blaschke product.
Let $A:=B \circ C$, let
$T$ be a unit disk automorphism and consider $T \circ B$.
Denote by $h : \D \to \R$ the least harmonic majorant of the subharmonic function $\log
|T \circ B|$. Note $h \le 0$.
Then $h \circ C$ is a harmonic majorant of $\log |T \circ
A|$ and $h \circ C\le 0$. Since $A$ is indestructible, $T \circ A$ is a Blaschke product, so that
by Lemma \ref{lem:2}, (a) $\Rightarrow$ (c), the least harmonic majorant of $\log |T \circ
A|$ ist $0$. It follows that $h \circ C = 0$. Since $C$ is not constant, the
function $h$,  the least harmonic majorant of $\log |T \circ B|$, is $0$. Hence
Lemma \ref{lem:2}, (c) $\Rightarrow$ (a), implies that $T \circ B$ is a Blaschke product.
This shows that $B$ is an indestructible Blaschke product.

\vfill
\hspace{0.9cm}\begin{minipage}{8cm}
Daniela Kraus and Oliver Roth\\
Department of Mathematics\\
University of W\"urzburg\\
97074 W\"urzburg\\
Germany\\
dakraus@mathematik.uni-wuerzburg.de\\
roth@mathematik.uni-wuerzburg.de

\end{minipage}

\end{document}